\documentclass[12pt,leqno]{amsart}
\usepackage{latexsym}
\usepackage{amsmath}
\usepackage{amsthm}
\usepackage{amssymb}
\usepackage{amsfonts}
\begin{document}
\theoremstyle{plain}
\newtheorem{MainThm}{Theorem}
\newtheorem{thm}{Theorem}[section]
\newtheorem{clry}[thm]{Corollary}
\newtheorem{prop}[thm]{Proposition}
\newtheorem{lemma}[thm]{Lemma}
\newtheorem{deft}[thm]{Definition}
\newtheorem{hyp}{Assumption}
\newtheorem*{conjecture}{Conjecture}

\newtheorem{claim}{Claim}[section]

\theoremstyle{definition}
\newtheorem{rem}[thm]{Remark}
\newtheorem*{acknow}{Acknowledgements}
\numberwithin{equation}{section}
\newcommand{\eps}{\varepsilon}
\renewcommand{\phi}{\varphi}
\renewcommand{\d}{\partial}
\newcommand{\re}{\mathop{\rm Re} }
\newcommand{\im}{\mathop{\rm Im}}
\newcommand{\mR}{\mathbb{R}}
\newcommand{\mC}{\mathbb{C}}
\newcommand{\mN}{\mathbb{N}} 
\newcommand{\mZ}{\mathbb{Z}} 
\newcommand{\mK}{\mathbb{K}}
\newcommand{\supp}{\mathop{\rm supp}}
\newcommand{\abs}[1]{\lvert #1 \rvert}
\newcommand{\norm}[1]{\lVert #1 \rVert}
\newcommand{\csubset}{\Subset}
\newcommand{\detg}{\lvert g \rvert}
\newcommand{\msetminus}{\setminus}

\newcommand{\mdiv}{\mathrm{div}}

\title[An inverse problem for the $p$-Laplacian]{An inverse problem for the $p$-Laplacian: boundary determination}

\author{Mikko Salo}
\address{Department of Mathematics and Statistics, University of Helsinki}
\email{mikko.salo@helsinki.fi}

\author{Xiao Zhong}
\address{Department of Mathematics, University of Jyv\"askyl\"a}
\email{xiao.x.zhong@jyu.fi}

\date{June 21, 2011}

\begin{abstract}
We study an inverse problem for nonlinear elliptic equations modelled after the $p$-Laplacian. It is proved that the boundary values of a conductivity coefficient are uniquely determined from boundary measurements given by a nonlinear Dirichlet-to-Neumann map. The result is constructive and local, and gives a method for determining the coefficient at a boundary point from measurements in a small neighborhood. The proofs work with the nonlinear equation directly instead of being based on linearization. In the complex valued case we employ complex geometrical optics type solutions based on $p$-harmonic exponentials, while for the real case we use $p$-harmonic functions first introduced by Wolff.
\end{abstract}

\maketitle

\section{Introduction} \label{section_introduction}

This article concerns inverse boundary value problems for nonlinear elliptic equations. In the case where the underlying equation is linear, a standard example is the inverse problem of Calder\'on \cite{C}. This problem is related to Electrical Impedance Tomography, a method proposed for medical and industrial imaging, where the objective is to determine the electrical conductivity of a medium by making voltage to current measurements on its boundary. There is an extensive theory concerning inverse boundary value problems for linear elliptic equations. We refer to \cite{U_IP} for a recent survey.

We recall the mathematical statement of the Calder\'on problem. Let $\Omega \subseteq \mR^n$, $n \geq 2$, be a bounded open set with $C^1$ boundary. We consider $\Omega$ as a medium which conducts electricity, with conductivity given by a positive function $\gamma \in L^{\infty}(\Omega)$. Assuming that there are no sources or sinks of current in $\Omega$, a voltage $f$ on the boundary induces a potential $u$ in the domain which (by Ohm's law) solves the Dirichlet problem 
\begin{equation*}
\left\{ \begin{array}{rll}
\mdiv(\gamma(x) \nabla u) &\!\!\!= 0 & \quad \text{in } \Omega, \\
u &\!\!\!= f & \quad \text{on } \partial \Omega.
\end{array} \right.
\end{equation*}
The boundary measurement corresponding to $f$, denoted by $\Lambda_{\gamma} f$, is the current at the boundary given by 
$$
\Lambda_{\gamma} f = \gamma \frac{\partial u}{\partial \nu}\Big|_{\partial \Omega}.
$$

We assume that one can prescribe many different voltages $f$ on the boundary, and then measure the corresponding boundary currents $\Lambda_{\gamma} f$. The map $\Lambda_{\gamma}$ is called the Dirichlet-to-Neumann map (DN map for short). Using a suitable weak formulation, $\Lambda_{\gamma}$ becomes a bounded linear operator from $W^{1/2,2}(\partial \Omega)$ to $W^{-1/2,2}(\partial\Omega)$ (where $W^{s,p}$ denotes the $L^p$ based Sobolev space with smoothness index $s$). The Calder\'on problem asks to determine the conductivity function $\gamma$ from knowledge of the operator $\Lambda_{\gamma}$.

There are several aspects of the Calder\'on problem that have been studied. We mention the following particular questions:

\begin{enumerate}
\item[1.] 
{\bf Boundary uniqueness.} If $\Lambda_{\gamma_1} = \Lambda_{\gamma_2}$, then $\gamma_1|_{\partial \Omega} = \gamma_2|_{\partial \Omega}$.
\item[2.]
{\bf Interior uniqueness.} If $\Lambda_{\gamma_1} = \Lambda_{\gamma_2}$, then $\gamma_1 = \gamma_2$ in $\Omega$.
\item[3.]
{\bf Stability.} If $\Lambda_{\gamma_1}$ and $\Lambda_{\gamma_2}$ are close, then $\gamma_1$ and $\gamma_2$ are close.
\item[4.]
{\bf Reconstruction.} Algorithm for determining $\gamma$ from $\Lambda_{\gamma}$.
\item[5.]
{\bf Partial data.} If $\Lambda_{\gamma_1}|_{\Gamma} = \Lambda_{\gamma_2}|_{\Gamma}$ for some $\Gamma \subseteq \partial \Omega$, then $\gamma_1 = \gamma_2$.
\end{enumerate}

As a rule, rather precise results for the above questions are available in the case $n=2$ \cite{Br}, \cite{AP}, \cite{BFR}, \cite{AP_boundary}, \cite{IUY}. Also in dimensions $n \geq 3$ many results have been obtained \cite{Br}, \cite{SU}, \cite{A}, \cite{N}, \cite{KSU}, however sharp conditions such as optimal regularity of the conductivity are in general not known. We refer again to the survey \cite{U_IP} for more details.

In contrast with the linear case, less is known about variants of the Calder\'on problem for nonlinear elliptic equations. In this paper we consider a particular nonlinear model based on the \emph{$p$-Laplace operator} 
$$
\Delta_p u = \mdiv(\abs{\nabla u}^{p-2} \nabla u), \quad 1 < p < \infty.
$$
The corresponding $p$-Laplace equation $\Delta_p u = 0$, whose solutions are called \emph{$p$-harmonic functions}, is a prototypical nonlinear equation in divergence form. It arises as the Euler-Lagrange equation for minimizing the $p$-Dirichlet energy 
$$
E_p(u) = \int_{\Omega} \abs{\nabla u}^p \,dx
$$
over all $u \in W^{1,p}(\Omega)$ with fixed boundary values. For more details on $p$-Laplace type operators we refer to the book \cite{HKM} and lecture notes \cite{dOI}, \cite{L}. Applications in fluid mechanics, plastic moulding, and image processing are discussed in \cite{Aronsson}, \cite{AronssonJanfalk}, \cite{Kuijper}.

Given a bounded open set $\Omega \subseteq \mR^n$, $n \geq 2$, having $C^1$ boundary, and for a positive function $\gamma \in L^{\infty}(\Omega)$, we consider the Dirichlet problem  
\begin{equation*}
\left\{ \begin{array}{rll}
\mdiv(\gamma(x) \abs{\nabla u}^{p-2} \nabla u) &\!\!\!= 0 & \quad \text{in } \Omega, \\
u &\!\!\!= f & \quad \text{on } \partial \Omega.
\end{array} \right.
\end{equation*}
This is also called the weighted $p$-Laplace equation \cite{HKM}. The nonlinear DN map is formally defined by 
\begin{equation*}
\Lambda_{\gamma}: f \mapsto \gamma(x) \abs{\nabla u}^{p-2} \nabla u \cdot \nu|_{\partial \Omega}
\end{equation*}
where $u$ is the unique solution with boundary values $f$ and $\nu$ is the outer unit normal to $\partial \Omega$. The precise definition of the DN map is given in Appendix \ref{section_appendix}.

Since the equation is nonlinear, one needs to make a distinction between real and complex valued solutions. We denote by $\Lambda_{\gamma}^{\mR}$ and $\Lambda_{\gamma}^{\mC}$ the DN maps acting on real and complex boundary values, respectively. Our main theorem is the following boundary uniqueness result.

\begin{thm} \label{thm_main1}
Let $\Omega \subseteq \mR^n$, $n \geq 2$, be a bounded open set having $C^1$ boundary, and let $\gamma_1, \gamma_2$ be positive continuous functions on $\overline{\Omega}$. If $\Lambda_{\gamma_1}^{\mR} = \Lambda_{\gamma_2}^{\mR}$, then $\gamma_1|_{\partial \Omega} = \gamma_2|_{\partial \Omega}$.
\end{thm}

The proof is constructive and local in the sense that we construct a sequence of explicit functions $f_M$ on the boundary, supported in an arbitrarily small neighborhood of a boundary point $x_0$, such that 
$$
\lim_{M \to \infty} \int_{\partial \Omega} \Lambda_{\gamma}(f_M) \bar{f}_M \,dS = \gamma(x_0).
$$
In fact, this result is true if $\gamma$ is in $L^{\infty}(\Omega)$ and continuous near $x_0$. The map $\Lambda_{\gamma}^{\mC}$ determines $\Lambda_{\gamma}^{\mR}$ trivially, and we obtain the following consequence for the complex valued case.

\begin{thm} \label{thm_main2}
Let $\Omega \subseteq \mR^n$, $n \geq 2$, be a bounded open set having $C^1$ boundary, and let $\gamma_1, \gamma_2$ be positive continuous functions on $\overline{\Omega}$. If $\Lambda_{\gamma_1}^{\mC} = \Lambda_{\gamma_2}^{\mC}$, then $\gamma_1|_{\partial \Omega} = \gamma_2|_{\partial \Omega}$.
\end{thm}

The reason for stating the theorems separately is that we actually first prove Theorem \ref{thm_main2} by employing solutions based on $p$-harmonic complex exponentials, and then establish Theorem \ref{thm_main1} by isolating the properties of complex exponentials needed for the proof and by making use of suitable real valued $p$-harmonic functions.

There are many works concerning boundary determination in the linear case where $p=2$. It was proved in \cite{KV} that the Taylor series at a boundary point of a smooth conductivity in a smooth domain is determined by the DN map. Another proof, based on pseudodifferential calculus and valid in many situations, was given in \cite{SU2}. For nonsmooth domains and conductivities there are boundary uniqueness results based on singular solutions \cite{A_boundary}, \cite{AG2} and direct methods involving explicit oscillating boundary values \cite{Br}, \cite{KY}, \cite{N_2D}, \cite{NT}.

Let us describe earlier results on nonlinear variants of the Calder\'on problem. One can consider the Dirichlet problem involving a more general nonlinear conductivity $a = a(x,z,q)$, 
\begin{equation*}
\left\{ \begin{array}{rll}
\mdiv(a(x,u,\nabla u) \nabla u) &\!\!\!= 0 & \quad \text{in } \Omega, \\
u &\!\!\!= f & \quad \text{on } \partial \Omega.
\end{array} \right.
\end{equation*}
Examples include 
\begin{enumerate}
\item[1.] 
(linear case) $a(x,z,q) = \gamma(x)$ for a positive $\gamma \in C^2(\overline{\Omega})$,
\item[2.]
(nonlinearity depending on $x$ and $u$) $a(x,z,q) = \gamma(x,z)$ where $\gamma \in C^2(\overline{\Omega} \times \mR)$ is positive,
\item[3.]
(derivative nonlinearity) $a(x,z,q) = \gamma(x,q)$ for suitable $\gamma$.
\end{enumerate}
In all these cases, the equation is uniformly elliptic and there is a unique solution in a suitable Sobolev space for any suitable boundary value $f$. We can then define the nonlinear DN map formally by 
\begin{equation*}
\Lambda_{\gamma}: f \mapsto a(x,u,\nabla u) \nabla u \cdot \nu|_{\partial \Omega}.
\end{equation*}

As explained above, many results are known for the linear case. Also the case where the nonlinearity depends on $x$ and $u$ is well understood. In fact, it was shown by Sun \cite{Sun_nonlinear} that this case reduces to the linear theory (the linearization idea is due to \cite{I1}). If $1 < p < \infty$, if $z \in \mR$, and if $\gamma^z(x) = \gamma(x,z)$, then one can prove that for any $f \in C^{2,\alpha}(\partial \Omega)$ 
$$
\lim_{t \to 0} \frac{\Lambda_{\gamma}(z+tf) - \Lambda_{\gamma}(z)}{t} = \Lambda_{\gamma^z}(f)
$$
in $W^{1-1/p,p}(\partial \Omega)$. As a consequence, if $\Lambda_{\gamma_1} = \Lambda_{\gamma_2}$ for two conductivities $\gamma_j = \gamma_j(x,u)$, then $\Lambda_{\gamma_1^z} = \Lambda_{\gamma_2^z}$ for all $z \in \mR$ and the interior uniqueness result in the linear case shows that $\gamma_1 = \gamma_2$ everywhere. Related results for other equations with nonlinearity depending on $x$ and $u$ appear in \cite{I1}, \cite{I2}, \cite{IN}, \cite{IS}, \cite{Sun2010}, \cite{SuU}.

For derivative nonlinearities it is possible to obtain some information from linearizations, see \cite{HS}, \cite{KN}, \cite{Sun_survey} for conductivity type equations and \cite{I3}, \cite{LW}, \cite{Sun2004} for related equations. These results are still based on linearizing the nonlinear DN map and applying the known uniqueness results for the linear case, and they apply to derivative nonlinearities of certain form. To deal with stronger nonlinearities one could hope for a method which works with the nonlinear equation directly. In this paper we introduce such a method, at least for the purposes of proving boundary determination results for $p$-Laplace type equations.

In the linear case most uniqueness results are based on special solutions of the conductivity equation, so called \emph{complex geometrical optics} solutions, which look like harmonic exponentials $e^{\rho \cdot x}$ where $\rho \in \mC^n$ and $\rho \cdot \rho = 0$. The first important observation is that also for the $p$-Laplacian there exist special complex solutions of the form $e^{\rho \cdot x}$ (see Lemma \ref{lemma_plaplace_approximate}). We prove Theorem \ref{thm_main2} by showing that one can perturb these $p$-harmonic complex exponentials to become solutions of the equation involving $\gamma$, concentrating near a boundary point. The proof is similar to the arguments of \cite{Br}, \cite{BrS} in the linear case and is actually not that difficult, making use of basic facts such as wellposedness for the Dirichlet problem, inequalities for $p$th powers of vectors, and Hardy's inequality. We then show Theorem \ref{thm_main1} by replacing the $p$-harmonic complex exponentials with certain real valued $p$-harmonic functions, introduced by Wolff \cite{W}, having similar properties as exponentials which allow the proof to go through.

It is an interesting question whether one can make progress on other aspects of inverse problems for $p$-Laplace type equations (such as interior uniqueness, stability, reconstruction, partial data) besides boundary uniqueness. Our results seem to suggest the possibility of a complex geometrical optics construction based on $p$-harmonic complex exponentials, or a corresponding construction in the real case using the functions of Wolff. One also expects the cases $n=2$ and $n \geq 3$ to be different. (In particular, when $n=2$ the $p$-harmonic equation is related to quasiregular mappings \cite{AIM} and also unique continuation for $p$-Laplace is known when $n=2$ \cite{BI} but it is not known for $n \geq 3$.) Another interesting direction would be to study more general equations modelled after the $p$-Laplacian, or to see if methods of this type are available for other physically relevant nonlinear equations. This paper is mainly intended to highlight a particular strongly nonlinear model for which Calder\'on type problems can be studied, and to give a first result in this setting.

This paper is organized as follows. Section \ref{section_introduction} is the introduction. The complex valued case is considered in Section \ref{section_complex} where Theorem \ref{thm_main2} is proved, while Section \ref{section_real} discusses the real valued case and proves Theorem \ref{thm_main1}. For the sake of completeness, there is an appendix containing standard material on inequalities for $p$th powers of vectors and on wellposedness and the Dirichlet-to-Neumann map for the equations considered in this paper. In the appendix we also make the observation that linearization of the DN map at constants does not give useful information of the conductivity, thus partly justifying the nonlinear methods used to prove the theorems.

\subsection*{Notation}

If $z, w \in \mC^n$ we write $z \cdot w = \sum_{j=1}^n z_j w_j$ for the dot product and $\abs{z} = (z \cdot \bar{z})^{1/2} = (\abs{\re(z)}^2 + \abs{\im(z)}^2)^{1/2}$ for the norm. If $1 \leq p \leq \infty$ the conjugate exponent is denoted by $p'$, so that $1/p + 1/p' = 1$. We write $W^{s,p}$ for the standard Sobolev spaces and $W^{1,p}_0$ for the closure of smooth compactly supported functions in $W^{1,p}$. The notation $A \lesssim B$ means that $A \leq CB$ for some constant $C > 0$ which is independent of asymptotic parameters (it typically only depends on $n$, $p$, and some choices of test functions). Similarly, $A \sim B$ means that $\frac{1}{C} A \leq B \leq C A$ for some constant $C > 0$.

\subsection*{Acknowledgements}

M.S.~is supported in part by the Academy of Finland, and is grateful to the Department of Mathematics and Statistics of the University of Jyv\"askyl\"a where part of this work was carried out. X.Z.~is supported by the Academy of Finland.

\section{Complex case} \label{section_complex}

We now prove Theorem \ref{thm_main2} concerning boundary determination from the DN map with complex boundary values. To convey the main idea without unnecessary complications, we will consider the case where $\Omega$ is a bounded open set in $\mR^n$ with $C^1$ boundary and $x_0$ is a point in $\partial \Omega$ such that $\partial \Omega$ is flat near $x_0$. (The non-flat case is covered when proving Theorem \ref{thm_main1} in Section \ref{section_real}.) By a translation and rotation, we may assume that $x_0 = 0$ and $\Omega \cap B(0,r) = \{ x \in B(0,r) \,;\, x_n > 0 \}$ for some small $r > 0$.

If $\gamma = \gamma(x) \in C(\overline{\Omega})$ is a positive function, consider the Dirichlet problem 
\begin{equation*}
\left\{ \begin{array}{rll}
\mdiv(\gamma(x) \abs{\nabla u}^{p-2} \nabla u) &\!\!\!= 0 & \quad \text{in } \Omega, \\
u &\!\!\!= f & \quad \text{on } \partial \Omega.
\end{array} \right.
\end{equation*}
The nonlinear DN map, acting on complex valued functions, is defined in the weak sense (see Appendix \ref{section_appendix}) by 
\begin{equation*}
\int_{\partial \Omega} \Lambda_{\gamma}(f) \bar{g} \,dS = \int_{\Omega} \gamma \abs{\nabla u}^{p-2} \nabla u \cdot \nabla \bar{v} \,dx, \quad f, g \in W^{1,p}(\Omega)
\end{equation*}
where $v$ is any function in $W^{1,p}(\Omega)$ with $v - g \in W^{1,p}_0(\Omega)$.

The main idea in the proof is to use a special solution to the nonlinear equation with coefficient $\gamma$ frozen at the boundary point $0$. This equation is just the $p$-Laplace equation, and the special solution is the following $p$-harmonic complex exponential.

\begin{lemma} \label{lemma_plaplace_approximate}
Let $h(x) = e^{\rho \cdot x}$ where $\rho = \alpha + i \beta$ with $\alpha, \beta \in \mR^n$. Then $h$ satisfies $\mdiv(\abs{\nabla h}^{p-2} \nabla h) = 0$ iff $(p-1)\abs{\alpha}^2 = \abs{\beta}^2$ and $\alpha \cdot \beta = 0$.
\end{lemma}
\begin{proof}
Since $\nabla h = \rho e^{\rho \cdot x}$, we have 
\begin{align*}
\mdiv(\abs{\nabla h}^{p-2} \nabla h) &= \mdiv(\abs{\rho}^{p-2} e^{(p-2)\alpha \cdot x} \rho e^{\rho \cdot x}) \\
 &= \mdiv(\abs{\rho}^{p-2} \rho e^{(p-1)\alpha \cdot x + i\beta \cdot x}) \\
 &= \abs{\rho}^{p-2} \rho \cdot ((p-1)\alpha + i\beta) e^{(p-1)\alpha \cdot x + i\beta \cdot x}.
\end{align*}
Here $\rho \cdot ((p-1)\alpha + i\beta) = (p-1)\abs{\alpha}^2 - \abs{\beta}^2 + ip \alpha \cdot \beta$, which proves the result.
\end{proof}

We wish to convert the $p$-harmonic function $e^{\rho \cdot x}$ into an exact solution of $\mdiv(\gamma \abs{\nabla u}^{p-2} \nabla u) = 0$ in $\Omega$ which concentrates near the boundary point $0$. To this end, define the function 
\begin{equation} \label{u0_definition}
u_0(x) = \eta_M(x) h_N(x)
\end{equation}
where $\eta_M(x) = \eta(Mx)$, $h_N(x) = h(Nx)$ where $M$ and $N$ are large positive numbers, $\eta \in C^{\infty}_c(\mR^n)$ is a nonnegative cutoff function with $\eta = 1$ for $\abs{x} \leq 1/2$ and $\eta = 0$ for $\abs{x} \geq 1$, and 
$$
h(x) = e^{(i\beta - e_n) \cdot x}
$$
with $\beta \in \mR^n$ satisfying $\abs{\beta}^2 = p-1$ and $\beta \cdot e_n = 0$. We will choose $N = N(M)$ so that $M/N \to 0$ as $M \to \infty$. The idea is that with these choices, since $h_N$ solves the equation with $\gamma$ frozen at $0$ and since $u_0$ is supported in the ball $B(0,1/M)$, $u_0$ becomes an approximate solution to the nonlinear equation in $\Omega$ when $M$ is large. Lemma \ref{lemma_u1_estimates} below gives a precise meaning to this statement.

We obtain an exact solution $u$ by solving the Dirichlet problem with boundary values $u_0$, 
\begin{equation*}
\left\{ \begin{array}{rll}
\mdiv(\gamma(x) \abs{\nabla u}^{p-2} \nabla u) &\!\!\!= 0 & \quad \text{in } \Omega, \\
u &\!\!\!= u_0 & \quad \text{on } \partial \Omega.
\end{array} \right.
\end{equation*}
Let $f = u_0|_{\partial \Omega}$. Then we have 
\begin{equation*}
\int_{\partial \Omega} \Lambda_{\gamma} (f) \bar{f} \,dS = \int_{\Omega} \gamma \abs{\nabla u}^{p-2} \nabla u \cdot \nabla \bar{u}_0 \,dx.
\end{equation*}
We write this as 
\begin{multline} \label{complex_identity_2}
\int_{\partial \Omega} \Lambda_{\gamma} (f) \bar{f} \,dS = \int_{\Omega} \gamma \abs{\nabla u_0}^p \,dx \\
 + \int_{\Omega} \gamma (\abs{\nabla u}^{p-2} \nabla u - \abs{\nabla u_0}^{p-2} \nabla u_0) \cdot \nabla \bar{u}_0 \,dx.
\end{multline}
Note that since $f$ is an explicit function, the left hand side is determined by the nonlinear DN map. We will recover the value of $\gamma$ at $0$ by taking the limit of this identity as $M \to \infty$. To analyze the limit, we need a simple lemma.

\begin{lemma} \label{lemma_zeta_estimates}
Let $\zeta \in C^{\infty}_c(B(0,1))$ and let $a \geq 0$. Then as $M \to \infty$ 
\begin{gather*}
M^{n-1} N \int_{\Omega} \zeta(Mx) e^{-p N x_n} \,dx \to \frac{1}{p} \int_{\mR^{n-1}} \zeta(x',0) \,dx', \\
\int_{\Omega} x_n^a \zeta(Mx) e^{-p N x_n} \,dx = O(M^{1-n} N^{-1-a}).
\end{gather*}
\end{lemma}
\begin{proof}
Follows from a direct computation.
\end{proof}

We compute the limit of the first term on the right hand side of \eqref{complex_identity_2}.

\begin{lemma} \label{lemma_u0_estimates}
We have as $M \to \infty$ 
\begin{equation*}
M^{n-1} N^{1-p} \int_{\Omega} \gamma \abs{\nabla u_0}^p \,dx \to c_p \gamma(0)
\end{equation*}
where $c_p = p^{\frac{p-2}{2}} \int_{\mR^{n-1}} \eta(x',0)^p \,dx'$. We also have
\begin{gather*}
\int_{\Omega} \abs{\nabla u_0}^p \,dx = O(M^{1-n} N^{p-1}), \\
\int_{\Omega} \abs{\eta_M \nabla h_N}^p \,dx = O(M^{1-n} N^{p-1}), \\
\int_{\Omega} \abs{\nabla u_0 - \eta_M \nabla h_N}^p \,dx = O(M^{1-n} N^{p-1}(M/N)^p).
\end{gather*}
\end{lemma}
\begin{proof}
Since $u_0 = \eta_M h_N$, we compute 
\begin{equation*}
\nabla u_0 = M \nabla \eta(M\,\cdot\,) h_N + \eta_M \nabla h_N.
\end{equation*}
Since $\nabla h_N = N(i\beta - e_n) e^{N(i\beta - e_n) \cdot x}=N(i\beta-e_n)h_N$, we have by Lemma \ref{lemma_zeta_estimates} 
\begin{gather*}
\norm{M \nabla \eta(M\,\cdot\,) h_N}_{L^p(\Omega)}^p = O(M^{1-n}N^{-1}M^p), \\
\norm{\eta_M \nabla h_N}_{L^p(\Omega)}^p = O(M^{1-n}N^{-1}N^p).
\end{gather*}
This shows the last three estimates since $M/N = o(1)$ as $M \to \infty$.

For the first statement, we use the inequality \eqref{z_convexity2} to conclude that 
\begin{align*}
 &\left\lvert \int_{\Omega} (\abs{\nabla u_0}^p - \abs{N(i\beta-e_n) \eta_M h_N}^p) \,dx \right\rvert \\
 &\leq p \int_{\Omega} \abs{M \nabla \eta(M\,\cdot\,) h_N}(\abs{\nabla u_0}^{p-1} + \abs{N(i\beta-e_n) \eta_M h_N}^{p-1}) \,dx \\
 &\leq p \norm{M \nabla \eta(M\,\cdot\,) h_N}_{L^p(\Omega)}(\norm{\nabla u_0}_{L^p(\Omega)}^{p-1} + \norm{N(i\beta-e_n) \eta_M h_N}_{L^p(\Omega)}^{p-1}) \\
 &= O(M^{1-n}N^{p-1}(M/N)).
\end{align*}
Using that $\abs{i\beta-e_n}^2 = p$, we have by Lemma \ref{lemma_zeta_estimates} 
\begin{align*}
\lim_{M\to\infty} M^{n-1} N^{1-p} \int_{\Omega} \abs{\nabla u_0}^p \,dx &= \lim_{M\to\infty} M^{n-1} N p^{p/2} \int_{\Omega} \abs{\eta_M h_N}^p \,dx \\
 &= p^{\frac{p-2}{2}} \int_{\mR^{n-1}} \eta(x',0)^p \,dx'.
\end{align*}
The result follows by writing $\gamma = \gamma(0) + (\gamma - \gamma(0))$ and by using the continuity of $\gamma$.
\end{proof}

We now move to the analysis of the second term on the right hand side of \eqref{complex_identity_2}. Writing $u = u_0 + u_1$, the next result shows that $\norm{\nabla u_1}_{L^p(\Omega)}$ is asymptotically smaller than $\norm{\nabla u_0}_{L^p(\Omega)}$. This may be interpreted so that $u_1$ is a small correction term which corrects the approximate solution $u_0$ into an exact solution $u$. The important facts for the proof are that $\Delta_p h = 0$ and that $u_0$ is supported near the boundary which makes it possible to use Hardy's inequality: if $\delta(x) = \text{dist}(x,\partial \Omega)$ then 
\begin{equation*}
\norm{v/\delta}_{L^p(\Omega)} \leq C \norm{\nabla v}_{L^p(\Omega)}, \quad v \in W^{1,p}_0(\Omega).
\end{equation*}

\begin{lemma} \label{lemma_u1_estimates}
As $M \to \infty$ 
\begin{equation*}
\int_{\Omega} \abs{\nabla u_1}^p \,dx = o(M^{1-n} N^{p-1}).
\end{equation*}
\end{lemma}
\begin{proof}
We will prove that
\begin{align}
I&=\int_\Omega (\abs{\nabla u}+\abs{\nabla u_0})^ {p-2}\abs{\nabla u_1}^2\, dx \notag \\
&\le o(M^{1-n}N^{p-1})+o(1)\int_{\Omega} \abs{\nabla u_1}^p \,dx. \label{add1}
\end{align}
%from which the lemma follows. 
To prove \eqref{add1}, we start with
$$
I \lesssim \int_\Omega \gamma (\abs{\nabla u}+\abs{\nabla u_0})^ {p-2}\abs{\nabla u_1}^2\, dx,
$$
since $\gamma$ is positive on $\overline{\Omega}$.  Then we invoke the inequality \eqref{zminuswp_largep}. Since $u_1 = u-u_0 \in W^{1,p}_0(\Omega)$ and since $u$ is a solution, we obtain that 
\begin{align*}
I &\lesssim \re \left[ \int_{\Omega} \gamma (\abs{\nabla u}^{p-2} \nabla u - \abs{\nabla u_0}^{p-2} \nabla u_0) \cdot (\nabla \bar{u} - \nabla \bar{u}_0) \,dx \right] \\
 &= - \re \left[ \int_{\Omega} \gamma \abs{\nabla u_0}^{p-2} \nabla u_0 \cdot \nabla \bar{u}_1 \,dx \right].
\end{align*}
The function $u_0$ is supported in the ball $B(0,1/M)$. Consequently, writing $\gamma = \gamma(0) + (\gamma - \gamma(0))$, we have 
\begin{align*}
I \lesssim &\left\lvert \int_{\Omega} \abs{\nabla u_0}^{p-2} \nabla u_0 \cdot \nabla \bar{u}_1 \,dx \right\rvert \\
 &+ \int_{B(0,1/M) \cap \Omega} \abs{\gamma - \gamma(0)} \abs{\nabla u_0}^{p-1} \abs{\nabla u_1} \,dx\\
 = & I_1+I_2.
\end{align*}
Integral $I_2$ is bounded by $\norm{\gamma - \gamma(0)}_{L^{\infty}(B(0,1/M) \cap \Omega)} \norm{\nabla u_0}_{L^p}^{p-1} \norm{\nabla u_1}_{L^p}$, which implies by Lemma \ref{lemma_u0_estimates}, the continuity of $\gamma$ and Young's inequality  that 
\begin{equation*}
I_2\le o(M^{1-n}N^{p-1})+o(1)\int_{\Omega} \abs{\nabla u_1}^p \,dx.
\end{equation*}
Then we estimate integral $I_1$ as follows.
At this point it is convenient to replace $\nabla u_0$ with $\eta_M \nabla h_N$ by writing 
\begin{multline*}
\int_{\Omega} \abs{\nabla u_0}^{p-2} \nabla u_0 \cdot \nabla \bar{u}_1 \,dx = \int_{\Omega} \abs{\eta_M \nabla h_N}^{p-2} \eta_M \nabla h_N \cdot \nabla \bar{u}_1 \,dx \\
 \quad + \int_{\Omega} (\abs{\nabla u_0}^{p-2} \nabla u_0 - \abs{\eta_M \nabla h_N}^{p-2} \eta_M \nabla h_N) \cdot \nabla \bar{u}_1 \,dx.
\end{multline*}
Integrating by parts, we obtain that 
\begin{multline*}
I_1 \lesssim \left\lvert \int_{\Omega} \mdiv(\abs{\eta_M \nabla h_N}^{p-2} \eta_M \nabla h_N ) \bar{u}_1 \,dx \right\rvert \\
 + \left\lvert \int_{\Omega} (\abs{\nabla u_0}^{p-2} \nabla u_0 - \abs{\eta_M \nabla h_N}^{p-2} \eta_M \nabla h_N) \cdot \nabla \bar{u}_1 \,dx \right\rvert.
\end{multline*}
In the first term on the right, we multiply and divide by $\delta$ (the distance to the boundary) and use the H\"older and Hardy inequalities so that 
\begin{multline*}
\left\lvert \int_{\Omega} \mdiv(\abs{\eta_M \nabla h_N}^{p-2} \eta_M \nabla h_N ) \bar{u}_1 \,dx \right\rvert \\
 \lesssim \norm{\delta \,\mdiv(\abs{\eta_M \nabla h_N}^{p-2} \eta_M \nabla h_N )}_{L^{p'}} \norm{\nabla u_1}_{L^p}.
\end{multline*}
The second term on the right can be estimated by \eqref{zpminus2z_difference}, and we have 
\begin{multline*}
\left\lvert \int_{\Omega} (\abs{\nabla u_0}^{p-2} \nabla u_0 - \abs{\eta_M \nabla h_N}^{p-2} \eta_M \nabla h_N) \cdot \nabla \bar{u}_1 \,dx \right\rvert \\
 \lesssim \int_{\Omega} (\abs{\nabla u_0} + \abs{\eta_M \nabla h_N})^{p-2} \abs{\nabla u_0 - \eta_M \nabla h_N} \abs{\nabla u_1} \,dx,
\end{multline*} 
which, by the H\"older inequality, is bounded by
\begin{multline*}
(\norm{\nabla u_0}_{L^p} + \norm{\eta_M \nabla h_N}_{L^p})^{p-2} \norm{\nabla u_0 - \eta_M \nabla h_N}_{L^p} \norm{\nabla u_1}_{L^p} \\
  = O((M^{1-n} N^{p-1})^{\frac{p-1}{p}} M/N) \norm{\nabla u_1}_{L^p},
\end{multline*}
when $p\ge 2$,  and by
\begin{multline*}
\norm{\nabla u_0 - \eta_M \nabla h_N}_{L^p}^{p-1} \norm{\nabla u_1}_{L^p} \\
=O((M^{1-n} N^{p-1})^{\frac{p-1}{p}} (M/N)^{p-1}) \norm{\nabla u_1}_{L^p},
\end{multline*}
when $1<p<2$.
In both cases, we used Lemma \ref{lemma_u0_estimates}. Since $M/N=o(1)$, we obtain that
\begin{multline*}
\left\lvert \int_{\Omega} (\abs{\nabla u_0}^{p-2} \nabla u_0 - \abs{\eta_M \nabla h_N}^{p-2} \eta_M \nabla h_N) \cdot \nabla \bar{u}_1 \right\rvert \\
\lesssim  o(M^{1-n}N^{p-1})+o(1)\int_{\Omega} \abs{\nabla u_1}^p \,dx.
\end{multline*}

Collecting these estimates together, we have proved that 
\begin{multline*}
I \lesssim \norm{\delta \,\mdiv(\abs{\eta_M \nabla h_N}^{p-2} \eta_M \nabla h_N )}_{L^{p'}}\norm{\nabla u_1}_{L^p}\\
+ o(M^{1-n} N^{p-1})+o(1)\int_{\Omega} \abs{\nabla u_1}^p \,dx.
\end{multline*}
We claim that as $M \to \infty$, 
\begin{equation}\label{add2}
\norm{\delta \,\mdiv(\abs{\eta_M \nabla h_N}^{p-2} \eta_M \nabla h_N )}_{L^{p'}}^{p'} = o(M^{1-n} N^{p-1}),
\end{equation}
from which estimate \eqref{add1} follows by Young's inequality. So, it remains to prove \eqref{add2}. 
Since $\eta_M$ and $h_N$ are explicit functions, this follows from a direct computation. Noting that $\mdiv(\abs{\nabla h_N}^{p-2} \nabla h_N) = \Delta_p h_N = N^p \Delta_p h = 0$, we have 
\begin{multline*}
\mdiv(\abs{\eta_M \nabla h_N}^{p-2} \eta_M \nabla h_N ) = \nabla(\eta_M^{p-1}) \cdot \abs{\nabla h_N}^{p-2} \nabla h_N \\
 = (p-1) \eta_M^{p-2} M \nabla \eta(M \,\cdot\,) N^{p-1} (\abs{\nabla h}^{p-2} \nabla h)(N \,\cdot\,).
\end{multline*}
Consequently, since $\delta(x) = x_n$, 
\begin{align*}
 &\norm{\delta \,\mdiv(\abs{\eta_M \nabla h_N}^{p-2} \eta_M \nabla h_N )}_{L^{p'}}^{p'} \\
 &\quad \lesssim M^{\frac{p}{p-1}} N^p \int_{B(0,1/M) \cap \Omega} x_n^{\frac{p}{p-1}} \abs{\nabla h(Nx)}^p \,dx \\
  &\quad \leq M^{\frac{p}{p-1}} N^{p-1-\frac{p}{p-1}} \int_0^{\infty} \int_{\abs{x'} \leq 1/M} x_n^{\frac{p}{p-1}} \abs{\nabla h(Nx',x_n)}^p \,dx \\
 &\quad \lesssim M^{\frac{p}{p-1}} N^{p-1-\frac{p}{p-1}} \int_0^{\infty} \int_{\abs{x'} \leq 1/M} x_n^{\frac{p}{p-1}} e^{-p x_n} \,dx' \,dx_n \\
 &= O(M^{\frac{p}{p-1} - n + 1} N^{p - 1 - \frac{p}{p-1}}).
\end{align*}
This is $O(M^{1-n} N^{p-1} (M/N)^{\frac{p}{p-1}}) = o(M^{1-n} N^{p-1})$ as required. This finishes the proof of \eqref{add2}, and hence that of \eqref{add1}. 
Now the lemma follows easily from \eqref{add1}. When $p\ge 2$, we have 
$$
I=\int_\Omega (\abs{\nabla u}+\abs{\nabla u_0})^ {p-2}\abs{\nabla u_1}^2\, dx\ge \int_\Omega\vert\nabla u_1\vert^p\, dx,
$$
which, together with \eqref{add1}, implies the desired estimate in the lemma. When $1<p<2$,
we have by H\"older's inequality
$$
\int_\Omega\vert \nabla u_1\vert^p\, dx\le I^{\frac{p}{2}}\left(\int_\Omega(\vert \nabla u\vert+\vert\nabla u_0\vert)^p\, dx\right)^{\frac{2-p}{2}},
$$
which implies the lemma.
\end{proof}

We now prove the following result, which immediately implies Theorem \ref{thm_main2} in the case where the boundary is flat near the point of interest.

\begin{prop}\label{propc}
If $\Omega$ is as above, there exists a sequence of explicit functions $(v_M) \subseteq C^{\infty}_c(\mR^n)$ such that their boundary values $f_M = v_M|_{\partial \Omega}$ satisfy $\supp(f_M) \subseteq B(0,1/M) \cap \partial \Omega$ and 
$$
\lim_{M \to \infty} \int_{\partial \Omega} \Lambda_{\gamma}(f_M) \bar{f}_M \,dS = \gamma(0).
$$
\end{prop}
\begin{proof}
If $f = u_0|_{\partial M}$ where $u_0$ is as in \eqref{u0_definition}, then \eqref{complex_identity_2} holds true. By Lemma \ref{lemma_u0_estimates}, we have 
$$
M^{n-1} N^{1-p} \int_{\Omega} \gamma \abs{\nabla u_0}^p \,dx \to c_p \gamma(0)
$$
where $c_p = p^{\frac{p-2}{2}} \int_{\mR^{n-1}} \eta(x',0)^p \,dx'$. By \eqref{zpminus2z_difference}, 
\begin{align*}
 &\left\lvert \int_{\Omega} \gamma (\abs{\nabla u}^{p-2} \nabla u - \abs{\nabla u_0}^{p-2} \nabla u_0) \cdot \nabla \bar{u}_0 \,dx \right\rvert \\
 &\quad \lesssim \int_{\Omega} (\abs{\nabla u} + \abs{\nabla u_0})^{p-2} \abs{\nabla u - \nabla u_0} \abs{\nabla u_0} \,dx \\
\intertext{If $p \geq 2$ then the H\"older inequality and Lemmas \ref{lemma_u0_estimates} and \ref{lemma_u1_estimates} imply that the last expression is bounded by}
 &\quad \lesssim (\norm{\nabla u}_{L^p} + \norm{\nabla u_0}_{L^p})^{p-2} \norm{\nabla u_1}_{L^p} \norm{\nabla u_0}_{L^p} \\
 &\quad = o(M^{1-n} N^{p-1}).
\end{align*}
If $1 < p < 2$ we obtain the same estimate from
\begin{align*}
 &\int_{\Omega} (\abs{\nabla u} + \abs{\nabla u_0})^{p-2} \abs{\nabla u - \nabla u_0} \abs{\nabla u_0} \,dx \\
 &\quad \leq \int_{\Omega} \abs{\nabla u - \nabla u_0}^{p-1} \abs{\nabla u_0} \,dx \\ 
 &\quad \lesssim \norm{\nabla u_1}_{L^p}^{p-1} \norm{\nabla u_0}_{L^p} \\
 &\quad = o(M^{1-n} N^{p-1}).
\end{align*}

Thus, if we define 
$$
v_M = \left( \frac{M^{n-1} N^{1-p}}{c_p} \right)^{1/2} u_0
$$
then the result follows.
\end{proof}

\section{Real valued case} \label{section_real}

Theorem \ref{thm_main1} concerning the nonlinear DN map acting on real valued functions is proved in this section. When proving Theorem \ref{thm_main2} we made use of $p$-harmonic complex exponentials. If $p=2$ there is little difference between real and complex solutions, since the real and imaginary parts of a complex solution can be used as real valued solutions. However, in the nonlinear case when $p \neq 2$ it is not clear how to obtain real valued solutions from complex ones. Therefore we cannot directly use the complex exponentials to establish Theorem \ref{thm_main1}, and the proof will require a real valued replacement for these functions.

Inspecting the proof of Theorem \ref{thm_main2} carefully, we see that the exact form of the exponential $h(x)$ is not essential. Rather, one needs that $h$ has certain properties, and it turns out that there exist real valued $p$-harmonic functions enjoying these properties. The following oscillatory solutions which decay exponentially in the $x_n$ variable were used by Wolff \cite[Section 3]{W}.

\begin{lemma} \label{lemma_functions_wolff}
The function $h(x) = e^{-x_n} a(x_1)$ satisfies $\Delta_p h = 0$ in $\mR^n_+$ if $a$ satisfies 
$$
a''(x_1) + V(a,a')a = 0,
$$
where 
$$
V(a,a') = \frac{(2p-3)(a')^2 + (p-1) a^2}{(p-1)(a')^2 + a^2}.
$$
Any solution $a$ is smooth and periodic with period $\lambda=\lambda(p)$, and satisfies $\int_0^\lambda a(x_1)\, dx_1=0$. 
\end{lemma}
Note that 
\begin{gather*}
\abs{h(x)} = e^{-x_n} \abs{a(x_1)}, \\
\abs{\nabla h(x)} = e^{-x_n} (a(x_1)^2 + a'(x_1)^2)^{1/2}.
\end{gather*}
We denote $h_N(x)=h(Nx)$ for $N>0$.

We begin the proof of Theorem \ref{thm_main1}. Assume that $\Omega \subseteq \mR^n$ is a bounded open set with $C^1$ boundary, and $x_0 \in \partial \Omega$. After a translation and rotation, we can assume that $x_0 = 0$ and the outer unit normal to $\partial \Omega$ at $0$ is the $n$th coordinate vector $-e_n$. We let $\rho$ be a boundary defining function for $\Omega$, that is, $\rho: \mR^n \to \mR$ is a $C^1$ function such that 
\begin{gather*}
\Omega = \{ x \in \mR^n \,;\, \rho(x) > 0 \}, \\
\partial \Omega = \{ x \in \mR^n \,;\, \rho(x) = 0 \}.
\end{gather*}
After scaling if necessary, we may assume that $\rho(0) = 0$ and $\nabla \rho(0) = e_n$. 
We define a map $f:\Omega\to \mR^n_+$ as 
$$
f(x)=(x^\prime, \rho(x)), \quad x=(x^\prime, x_n)\in \Omega.
$$
This map is invertible and close to the identity map in $B(0,1/M)\cap \Omega$ for $M$ large since $\rho$ is a $C^1$ function, which implies that as $M\to \infty$
\begin{equation}\label{add3}
\sup_{x\in B(0,1/M)\cap\Omega} \abs{\nabla \rho(x) - e_n} = o(1), \ \ \sup_{x\in B(0,1/M)\cap\Omega}\vert Df(x)-I\vert=o(1),
\end{equation}
where $I$ is the identity matrix.

Similarly as in Section \ref{section_complex}, we employ an approximate $p$-harmonic function 
$$
u_0(x) = \eta_M(x) \tilde{h}_N(x)
$$
where $\eta_M(x) = \eta(Mx)$ is the same cutoff function as in Section 2,  and 
$$
\tilde{h}_N(x) =h_N(f(x))= e^{-N\rho(x)} a(Nx_1),
$$
where $a$ is as in Lemma \ref{lemma_functions_wolff} and is not identically zero.

Note that $u_0 \in C^1_c(\mR^n)$. We choose a solution $u$ by requiring that 
\begin{equation*}
\left\{ \begin{array}{rll}
\mdiv(\gamma(x) \abs{\nabla u}^{p-2} \nabla u) &\!\!\!= 0 & \quad \text{in } \Omega, \\
u &\!\!\!= u_0 & \quad \text{on } \partial \Omega.
\end{array} \right.
\end{equation*}
Writing $f = u_0|_{\partial \Omega}$, we observe that the following real valued analogue of \eqref{complex_identity_2} holds:
\begin{multline} \label{real_identity_2}
\int_{\partial \Omega} \Lambda_{\gamma} (f) f \,dS = \int_{\Omega} \gamma \abs{\nabla u_0}^p \,dx \\
 + \int_{\Omega} \gamma (\abs{\nabla u}^{p-2} \nabla u - \abs{\nabla u_0}^{p-2} \nabla u_0) \cdot \nabla u_0 \,dx.
\end{multline}
As in Section 2, will recover the value of $\gamma$ at 0 by taking the limit of this identity as $M\to \infty$. To estimate the limits, we need the following lemmas,
Lemma \ref{lemma_real_zeta_estimates} and Lemma \ref{lemma_real_u0_estimates}, which are analogous to Lemma \ref{lemma_zeta_estimates} and Lemma \ref{lemma_u0_estimates}, respectively. 
As in Section 2, $M/N=o(1)$ as $M\to \infty$ and 
$\delta(x)=\text{dist}(x, \partial\Omega)$.

\begin{lemma} \label{lemma_real_zeta_estimates}
Let $\zeta \in C^{\infty}_c(B(0,2))$ and $\beta \geq 0$. 
Let $g:{\mathbb R}\to {\mathbb R}$ be a bounded periodic function with period $\lambda>0$. 
Then as $M \to \infty$ 
\begin{gather*}
M^{n-1} N \int_{\Omega} \zeta(Mx) e^{-p N \rho(x)} g(Nx_1) \,dx \to \frac{c}{p} \int_{\mR^{n-1}} \zeta(x',0) \,dx', \\
\int_{\Omega} \delta(x)^\beta \zeta(Mx) e^{-p N \rho(x)} \,dx = O(M^{1-n} N^{-1-\beta}),
\end{gather*}
where $c=\frac{1}{\lambda}\int_0^\lambda g(t)\, dt$.
\end{lemma}
\begin{proof}
We make the change of variables $x = f^{-1}(y)$, and \eqref{add3} together with easy computations shows that it is enough to prove these claims when $\Omega$ is flat near $0$. We only prove the first claim and compute 
\begin{align*}
 &M^{n-1} N \int_{\mR^n_+} \zeta(My) e^{-pNy_n} g(Ny_1) \,dy \\
 &\quad = \int_0^{\infty} \left( \int_{\abs{y'} \leq 2} \zeta(y',\frac{M}{N}y_n) g(\frac{N}{M} y_1) \,dy' \right) e^{-py_n} \,dy_n \\
 &\quad = \int_0^{\infty} \left( \int_{\abs{y'} \leq 2} (\zeta(y',0) + O(\frac{M}{N} y_n))g(\frac{N}{M} y_1) \,dy' \right) e^{-py_n} \,dy_n.
\end{align*}
Writing $R = N/M$ and $\psi(y_1) = \int_{\mR^{n-2}} \zeta(y_1,y'',0) \,dy''$, we have
\begin{align*}
 &\int_{\mR^{n-1}} \zeta(y',0) g(Ry_1) \,dy' = \frac{1}{R} \int_{-\infty}^{\infty} \psi(t/R)g(t) \,dt \\
 &\quad = \frac{1}{R} \sum_{k=-\infty}^{\infty} \int_{k\lambda}^{(k+1)\lambda} \psi(t/R) g(t) \,dt \\
 &\quad = \frac{1}{\lambda} \int_0^{\lambda} \left[ \sum_{k=-\infty}^{\infty} \psi((t+k\lambda)/R) \frac{\lambda}{R} \right] g(t) \,dt
\end{align*}
where the expression in brackets is a Riemann sum converging uniformly to $\int_{\mR^{n-1}} \zeta(y',0) \,dy'$ as $M \to \infty$. The result follows.
\end{proof}

\begin{lemma} \label{lemma_real_u0_estimates}
We have as $M \to \infty$ 
\begin{equation*}
M^{n-1} N^{1-p} \int_{\Omega} \gamma \abs{\nabla u_0}^p \,dx \to c_p \gamma(0)
\end{equation*}
where 
$$
c_p = \frac{K}{p}\int_{\mR^{n-1}} \eta(x',0)^p \,dx', \quad K=\frac{1}{\lambda}\int_0^\lambda\left(a(t)^2+a^\prime(t)^2\right)^{\frac{p}{2}}\, dt.
$$ 
We also have 
\begin{gather*}
\int_{\Omega} \abs{\nabla u_0}^p \,dx = O(M^{1-n} N^{p-1}), \\
\int_{\Omega} \abs{\eta_M \nabla \tilde{h}_N}^p \,dx = O(M^{1-n} N^{p-1}), \\
\int_{\Omega} \abs{\nabla u_0 - \eta_M \nabla \tilde{h}_N}^p \,dx = O(M^{1-n} N^{p-1}(M/N)^p).
\end{gather*}
\end{lemma}
\begin{proof}
Noting that $\nabla u_0 = M \nabla \eta(M \,\cdot\,) \tilde{h}_N + \eta_M \nabla \tilde{h}_N$, the last three estimates follow from Lemma \ref{lemma_real_zeta_estimates}. For the first claim we observe that 
$$
\nabla \tilde{h}_N(x) - \nabla h_N(f(x)) = -N e^{-N \rho(x)}(\nabla \rho(x) - \nabla \rho(0)) a(N x_1).
$$
Thus by the inequality \eqref{z_convexity2}, when computing the limit of $\int_{\Omega} \gamma \abs{\nabla u_0}^p \,dx$ as $M \to \infty$ we may replace $\gamma$ by $\gamma(0)$ and $\nabla u_0$ by $\eta_M(\,\cdot\,) \nabla h_N(f(\,\cdot\,))$. Lemma \ref{lemma_real_zeta_estimates} gives the required expression for the limit.
\end{proof}

Next we write $u = u_0 + u_1$ and show that $\nabla u_1$ is asymptotically smaller than $\nabla u_0$ in the $L^p$ norm.

\begin{lemma} \label{lemma_real_u1_estimates}
As $M \to \infty$
$$
\int_{\Omega} \abs{\nabla u_1}^p \,dx = o(M^{1-n} N^{p-1}).
$$
\end{lemma}
\begin{proof}
The proof is analogous to the proof of Lemma \ref{lemma_u1_estimates}. We start with
\begin{align*}
 I&=\int_{\Omega} (\abs{\nabla u}+\abs{\nabla u_0})^{p-2}\abs{\nabla u_1}^2 \,dx\\
 &\lesssim \int_{\Omega} \gamma (\abs{\nabla u}+\abs{\nabla u_0})^{p-2}\abs{\nabla u_1}^2 \,dx \\
 &\lesssim \int_{\Omega} \gamma (\abs{\nabla u}^{p-2} \nabla u - \abs{\nabla u_0}^{p-2} \nabla u_0) \cdot (\nabla u - \nabla u_0) \,dx \\
 & = - \int_{\Omega} \gamma \abs{\nabla u_0}^{p-2} \nabla u_0 \cdot \nabla u_1 \,dx \\
 &= -\gamma(0) \int_{\Omega} \abs{\nabla u_0}^{p-2} \nabla u_0 \cdot \nabla u_1 \,dx + o((M^{1-n}N^{p-1})^{\frac{p-1}{p}})\norm{\nabla u_1}_{L^p}.
\end{align*}
In the last step we used Lemma  \ref{lemma_real_u0_estimates}, the continuity of $\gamma$ and the fact that $u_0$ is supported in $B(0,1/M)$.

We continue by writing 
\begin{multline*}
\int_{\Omega} \abs{\nabla u_0}^{p-2} \nabla u_0 \cdot \nabla u_1 \,dx = \int_{\Omega} \abs{\eta_M \nabla \tilde h_N}^{p-2} \eta_M \nabla \tilde h_N \cdot \nabla u_1 \,dx \\
 + \int_{\Omega} (\abs{\nabla u_0}^{p-2} \nabla u_0 - \abs{\eta_M \nabla \tilde h_N}^{p-2} \eta_M \nabla \tilde h_N) \cdot \nabla u_1 \,dx.
\end{multline*}
For the last term, we can apply Lemma  \ref{lemma_real_u0_estimates} and proceed in the same way as in the proof of Lemma \ref{lemma_u1_estimates} to obtain that
\begin{align*}
 &\left\lvert \int_{\Omega} (\abs{\nabla u_0}^{p-2} \nabla u_0 - \abs{\eta_M \nabla \tilde  h_N}^{p-2} \eta_M \nabla \tilde h_N) \cdot \nabla u_1 \,dx \right\rvert \\
 &\qquad =o((M^{1-n}N^{p-1})^{\frac{p-1}{p}})\norm{\nabla u_1}_{L^p}.
\end{align*}
For the first term in the right hand,  we estimate in the following way. 
Note that $\tilde h_N(x)=h_N(f(x))$. We make a change of variables and use that $\nabla \tilde{h}_N(x)$ can be replaced by $\nabla h_N(f(x))$ up to small error to obtain 
\begin{align*}
 & \int_{\Omega} \abs{\eta_M \nabla \tilde h_N}^{p-2} \eta_M \nabla \tilde h_N \cdot \nabla u_1 \,dx \\
 &\quad =\int_{\mR^n_+} \eta_M^{p-1}(f^{-1}(y)) \abs{\nabla h_N}^{p-2} \nabla h_N\cdot\nabla u_1(f^{-1}(y)) \,dy \\
 &\qquad + o((M^{1-n} N^{p-1})^{\frac{p-1}{p}})\norm{\nabla u_1}_{L^p}.
\end{align*}
Here we also used \eqref{add3}. Then we can estimate the first term in the right hand side by Hardy's inequality in the same way as in the proof of Lemma \ref{lemma_u1_estimates}. Eventually we obtain the required estimate for $I$,
$$
I = o((M^{1-n} N^{p-1})^{\frac{p-1}{p}})\norm{\nabla u_1}_{L^p} = o(M^{1-n}N^{p-1})+o(1)\int_{\Omega} \abs{\nabla u_1}^p \,dx.
$$
We can conclude the proof as in Lemma \ref{lemma_u1_estimates}.
\end{proof}

Theorem \ref{thm_main1} now follows from the next result.

\begin{prop}\label{propr}
If $\Omega$ is as above, there exists a sequence of explicit real valued functions $(v_M) \subseteq C^1_c(\mR^n)$ such that their boundary values $f_M = v_M|_{\partial \Omega}$ satisfy $\supp(f_M) \subseteq B(0,1/M) \cap \partial \Omega$ and 
$$
\lim_{M \to \infty} \int_{\partial \Omega} \Lambda_{\gamma}(f_M) f_M \,dS = \gamma(0).
$$
\end{prop}
\begin{proof}
As in the proof of Proposition \ref{propc}, we set 
$$
v_M = \left(\frac{M^{n-1}N^{1-p}}{c_p}\right)^{\frac{1}{2}}u_0,
$$
where the constant $c_p$ is the same as in Lemma \ref{lemma_real_u0_estimates}. Then the result follows.
\end{proof}

\appendix

\section{Inequalities and wellposedness} \label{section_appendix}

We first record some standard inequalities for $p$th powers of vectors. The first set of inequalities is valid for $1 < p < \infty$ and for all $\xi, \eta \in \mR^n$.
\begin{gather}
\abs{\eta}^p \geq \abs{\xi}^p + p \abs{\xi}^{p-2} \xi \cdot (\eta-\xi), \\
\abs{\abs{\xi}^p - \abs{\eta}^p} \leq p (\abs{\xi}^{p-1} + \abs{\eta}^{p-1}) \abs{\xi-\eta}, \\
\abs{\abs{\xi}^{p-2} \xi - \abs{\eta}^{p-2} \eta} \lesssim (\abs{\xi} + \abs{\eta})^{p-2} \abs{\xi-\eta}, \\
(\abs{\xi}^{p-2} \xi - \abs{\eta}^{p-2} \eta) \cdot (\xi-\eta) \sim (\abs{\xi} + \abs{\eta})^{p-2} \abs{\xi-\eta}^2.
\end{gather}
These results easily imply similar inequalities for complex vectors upon separating real and imaginary parts. If $z, w \in \mC^n$ and $1 < p < \infty$ we have 
\begin{gather}
\abs{w}^p \geq \abs{z}^p + p \abs{z}^{p-2} \re[z \cdot (\bar{w} - \bar{z})], \label{z_convexity} \\ 
\abs{\abs{z}^p - \abs{w}^p} \leq p (\abs{z}^{p-1} + \abs{w}^{p-1}) \abs{z-w}, \label{z_convexity2} \\
\abs{\abs{z}^{p-2} z - \abs{w}^{p-2} w} \lesssim (\abs{z} + \abs{w})^{p-2} \abs{z-w}. \label{zpminus2z_difference} \\
(\vert z\vert+\vert w\vert)^{p-2}\abs{z-w}^2 \sim \re \left[ (\abs{z}^{p-2} z - \abs{w}^{p-2} w) \cdot (\bar{z} - \bar{w}) \right]. \label{zminuswp_largep}
\end{gather}
%while for $1 < p \leq 2$ we will use that 
%\begin{equation}
%\abs{z-w}^p \lesssim \left( \re \left[ (\abs{z}^{p-2} z - \abs{w}^{p-2} w) \cdot (\bar{z} - \bar{w}) \right] \right)^{p/2} (\abs{z}^p + \abs{w}^p)^{1-p/2}.
%\end{equation}

Assume that $\Omega \subseteq \mR^n$ is a bounded open set, and that $\gamma \in L^{\infty}(\Omega)$ is a positive function. Let $\mK$ be either $\mR$ or $\mC$. We consider the wellposedness of the Dirichlet problem for the $p$-Laplace equation, 
\begin{equation*}
\left\{ \begin{array}{rll}
\mdiv(\gamma(x) \abs{\nabla u}^{p-2} \nabla u) &\!\!\!= 0 & \quad \text{in } \Omega, \\
u &\!\!\!= f & \quad \text{on } \partial \Omega.
\end{array} \right.
\end{equation*}
We look for a weak solution $u \in W^{1,p}(\Omega ; \mK)$, so that 
$$
\int_{\Omega} \gamma \abs{\nabla u}^{p-2} \nabla u \cdot \nabla \varphi \,dx = 0 \quad \text{for all $\varphi \in W^{1,p}_0(\Omega ; \mK)$}.
$$
Here $f \in W^{1,p}(\Omega ; \mK)$, and the boundary condition is interpreted so that $u-f \in W^{1,p}_0(\Omega ; \mK)$. We think of $f$ belonging to the abstract trace space $W^{1,p}(\Omega ; \mK) / W^{1,p}_0(\Omega ; \mK)$, and also write $u|_{\partial \Omega} = f$. If $\Omega$ has sufficiently nice (say Lipschitz) boundary, the trace space can be identified with the Besov space \cite{JK} 
$$
W^{1,p}(\Omega ; \mK) / W^{1,p}_0(\Omega ; \mK) \approx B^{1-1/p}_{pp}(\partial \Omega ; \mK).
$$

\begin{prop}
Given any $f \in W^{1,p}(\Omega ; \mK)$ the above Dirichlet problem has a unique solution $u \in W^{1,p}(\Omega ; \mK)$ satisfying 
$$
\norm{u}_{W^{1,p}} \leq C \norm{f}_{W^{1,p}}
$$
with $C$ independent of $f$.
\end{prop}
\begin{proof}
The proof is completely standard and is usually given for the case $\mK = \mR$ \cite{L}. We show that the same proof works for $\mK = \mC$.

We will show the solution is obtained as the unique minimizer of the energy functional 
$$
E(v) = \int_{\Omega} \gamma \abs{\nabla v}^p \,dx
$$
among all $v$ in the set $\mathcal{A} = \{ v \in W^{1,p}(\Omega ; \mC) \,;\, v-f \in W^{1,p}_0(\Omega ; \mC) \}$. In fact, let $u$ be a minimizer. Then for any $\varphi \in W^{1,p}_0(\Omega ; \mC)$ 
\begin{align*}
E(u + t\varphi) &= \int_{\Omega} \gamma (\abs{\nabla u + t \nabla \varphi}^2)^{p/2} \,dx \\
 &= \int_{\Omega} \gamma (\abs{\nabla u}^2 + 2t \re(\nabla u \cdot \overline{\nabla \varphi}) + t^2 \abs{\nabla \varphi}^2)^{p/2} \,dx
\end{align*}
and 
\begin{align*}
\frac{d}{dt} E(u+t\varphi)\Big|_{t=0} = p \int_{\Omega} \gamma \abs{\nabla u}^{p-2} \re(\nabla u \cdot \overline{\nabla \varphi}) \,dx.
\end{align*}
Since $u$ is a minimizer the last expression is zero, and choosing $\varphi$ purely real or purely imaginary gives that $u$ is a solution, 
$$
\mdiv(\gamma \abs{\nabla u}^{p-2} \nabla u) = 0.
$$
Conversely, if $u$ is a solution with $u-f \in W^{1,p}_0(\Omega ; \mC)$ then for any $v \in \mathcal{A}$ we have by \eqref{z_convexity}
\begin{align*}
E(v) &= \int_{\Omega} \gamma \abs{\nabla v}^p \,dx \\
 &\geq \int_{\Omega} \gamma \left( \abs{\nabla u}^p + p \abs{\nabla u}^{p-2} \re\left[ \nabla u \cdot (\nabla \bar{v} - \nabla \bar{u}) \right] \right) \,dx.
\end{align*}
Since $v-u \in W^{1,p}_0(\Omega ; \mC)$ one has $E(v) \geq E(u)$ and $u$ is a minimizer.

It is thus enough to show that the energy functional has a unique minimizer in $\mathcal{A}$. We begin with uniqueness, and suppose that $u_1, u_2 \in \mathcal{A}$ are two distinct minimizers. Then $\nabla u_1 \neq \nabla u_2$ in a set of positive measure (since otherwise $u_1 - u_2 \in W^{1,p}_0(\Omega ; \mC)$ is constant on each component and thus zero by the boundary condition). The function $\frac{u_1+u_2}{2}$ also belongs to $\mathcal{A}$, and the strict convexity of $x \mapsto \abs{x}^p$ implies that 
$$
\left\lvert \frac{\nabla u_1 + \nabla u_2}{2} \right\rvert^p \leq \frac{\abs{\nabla u_1}^p + \abs{\nabla u_2}^p}{2}
$$
where the inequality is strict in a set of positive measure. We obtain 
\begin{align*}
E(u_1) &= \int_{\Omega} \gamma \abs{\nabla u_1}^p \,dx \leq \int_{\Omega} \gamma \left\lvert \frac{\nabla u_1 + \nabla u_2}{2} \right\rvert^p \,dx \\
 &< \int_{\Omega} \gamma \frac{\abs{\nabla u_1}^p + \abs{\nabla u_2}^p}{2} \,dx.
\end{align*}
This implies $E(u_1) < E(u_2) \leq E(u_1)$, which is a contradiction.

To show the existence of a minimizer, we let $(v_j) \subseteq \mathcal{A}$ be a sequence such that $E(v_j) \to E_0$ where $E_0 = \inf_{v \in \mathcal{A}} E(v)$. Since $f \in \mathcal{A}$, we have $E_0 < \infty$ and 
$$
\int_{\Omega} \abs{\nabla v_j} \,dx \leq C < \infty.
$$
Using that $v_j - f \in W^{1,p}_0(\Omega ; \mC)$, the Poincar\'e inequality implies 
$$
\norm{v_j - f}_{W^{1,p}} \leq C \norm{\nabla(v_j - f)}_{L^p} \leq C < \infty.
$$
By weak compactness there exists $u \in W^{1,p}(\Omega ; \mC)$ with $v_j \rightharpoonup u$ and $\nabla v_j \rightharpoonup \nabla u$ weakly in $L^p(\Omega ; \mC)$, and since $W^{1,p}_0$ is closed we have $u-f \in W^{1,p}_0(\Omega ; \mC)$. The function $u$ is a minimizer: by \eqref{z_convexity} 
\begin{align*}
E(v_j) &= \int_{\Omega} \gamma \abs{\nabla v_j}^p \,dx \\
 &\geq \int_{\Omega} \gamma \left( \abs{\nabla u}^p + p \abs{\nabla u}^{p-2} \re\left[ \nabla u \cdot (\nabla \bar{v_j} - \nabla \bar{u}) \right] \right) \,dx
\end{align*}
and by weak convergence $E_0 = \lim_{j \to \infty} E(v_j) \geq E(u)$. The norm bound for $u$ follows from the fact that $E(u) \leq E(f)$ and from the Poincar\'e inequality.
\end{proof}

Finally, we discuss a nonlinear Dirichlet-to-Neumann map (DN map) for the equation considered above. Write $X = W^{1,p}(\Omega)/W^{1,p}_0(\Omega)$ and $X'$ for the dual of $X$. If $f \in X$, the DN map is formally defined by 
$$
\Lambda_{\gamma}(f) = \gamma \abs{\nabla u_f}^{p-2} \nabla u_f \cdot \nu|_{\partial \Omega}
$$
where $u_f \in W^{1,p}(\Omega)$ is the unique solution of $\mdiv(\gamma \abs{\nabla u}^{p-2} \nabla u) = 0$ in $\Omega$ with $u-f \in W^{1,p}_0(\Omega)$, and $\nu$ is the outer unit normal vector of $\partial \Omega$. A formal integration by parts gives that 
$$
\langle \Lambda_{\gamma}(f), g \rangle_{\partial \Omega} = \int_{\partial \Omega} \Lambda_{\gamma}(f) g \,dS = \int_{\Omega} \gamma \abs{\nabla u_f}^{p-2} \nabla u_f \cdot \nabla g \,dx.
$$
The last identity can be taken as the weak definition of the DN map for $f \in X$. One has 
$$
\abs{\langle \Lambda_{\gamma}(f), g \rangle_{\partial \Omega}} \leq C \norm{\nabla u_f}_{L^p}^{p-1} \norm{\nabla g}_{L^p} \leq C \norm{f}_{W^{1,p}}^{p-1} \norm{g}_{W^{1,p}}.
$$
Consequently $\Lambda_{\gamma}$ maps $X$ to $X'$ and 
$$
\norm{\Lambda_{\gamma}(f)}_{X'} \leq C \norm{f}_X^{p-1}.
$$
%If $\Omega$ is a bounded domain with smooth boundary and $\gamma$ is smooth, it is known that there exists $\alpha > 0$ so that for any $f \in C^{1,\alpha}(\partial \Omega)$ the Dirichlet problem has a unique solution $u_f \in C^{1,\alpha}(\overline{\Omega})$ with 
%$$
%\norm{u_f}_{C^{1,\alpha}(\overline{\Omega})} \leq C \norm{f}_{C^{1,\alpha}(\partial \Omega)}.
%$$
%Consequently, for this $\alpha$ the DN map takes $C^{1,\alpha}(\partial \Omega)$ into $C^{\alpha}(\partial \Omega)$. (This result is optimal in the sense that there exist smooth boundary values for which the corresponding solutions are not better than $C^{1,\alpha}$.)

Some earlier results for Calder\'on type inverse problems for nonlinear equations have been based on studying the G{\^a}teaux derivatives of the nonlinear map $\Lambda_{\gamma}$ at boundary values $f = z$ where $z$ is a constant. However, for the equation above the solution with boundary values $z+tf$ (where $t > 0$) is $u_{z+tf} = z + tu_f$, and 
$$
\Lambda_{\gamma}(z+tf) = \gamma \abs{\nabla(z+tu_f)}^{p-2} \nabla(z+tu_f) \cdot \nu|_{\partial \Omega} = t^{p-1} \Lambda_{\gamma}(f).
$$
Consequently the G{\^a}teaux derivatives of $\Lambda_{\gamma}$ at $z$ do not even exist if $1 < p < 2$, and also when $p > 2$ all the higher order G{\^a}teaux derivatives which exist are either $0$ or equal to $\Lambda_{\gamma}$. Thus for the $p$-Laplace type equation considered in this paper, arguments based on linearizing the map $\Lambda_{\gamma}$ at constants do not yield any information and one needs other methods.


\begin{thebibliography}{99}

\bibitem{A} G.~Alessandrini, \textit{{Stable determination of conductivity by boundary measurements}}, Appl. Anal. \textbf{27} (1988), 153--172.
\bibitem{A_boundary} G.~Alessandrini, \textit{{Singular solutions of elliptic equations and the determination of conductivity by boundary measurements}}, J. Diff. Eq. \textbf{84} (1990), 252--272.
\bibitem{AG2} G.~Alessandrini, R.~Gaburro, \textit{{The local Calder\'on problem and the determination at the boundary of the conductivity}}, Comm. PDE \textbf{34} (2009), 918--936.
\bibitem{Aronsson} G.~Aronsson, \textit{{On p-harmonic functions, convex duality and an asymptotic formula for injection mould filling}}, European J. Appl. Math. \textbf{7} (1996), 417--437.
\bibitem{AronssonJanfalk} G.~Aronsson, U. Janfalk, \textit{{On Hele-Shaw flow of power-law fluids}}, European J. Appl. Math. \textbf{3} (1992), 343--366.
\bibitem{AIM} K.~Astala, T.~Iwaniec, and G.~Martin, \textit{{Elliptic partial differential equations and quasiconformal mappings in the plane}}, Princeton University Press, 2009.
\bibitem{AP_boundary} K.~Astala, L.~P\"aiv\"arinta, \textit{A boundary integral equation for Calder{\'o}n's inverse conductivity problem}, Collect. Math. (2006), Vol. Extra, 127--139.
\bibitem{AP} K.~Astala, L.~P\"aiv\"arinta, \textit{Calder{\'o}n's inverse conductivity problem in the plane}, Ann. of Math. \textbf{163} (2006), 265--299.
\bibitem{BFR} T.~Barcel\'o, D.~Faraco, and A.~Ruiz, \textit{{Stability of Calder{\'o}n inverse conductivity problem in the plane}}, J. Math. Pures Appl. \textbf{88} (2007), 522--556.
\bibitem{BI} B.~Bojarski and T.~Iwaniec, \textit{{$p$-harmonic equation and quasiregular mappings}}, Partial Differential Equations (Warsaw 1984), pp. 25-38, Banach Center Publications 19 (1987).
\bibitem{Br} R.~M.~Brown, \emph{{Recovering conductivity at the boundary from the Dirichlet to Neumann map: a pointwise result}}, J. Inverse Ill-Posed Probl. \textbf{9} (2001), 567--574.
\bibitem{BrS}
R.~M. Brown and M.~Salo, \emph{Identifiability at the boundary for first-order
  terms}, Appl. Anal. \textbf{85} (2006), no.~6-7, 735--749.
%\bibitem{Bu} A.~L.~Bukhgeim, \textit{Recovering a potential from Cauchy data in the two-dimensional case}, J. Inverse Ill-posed Probl., \textbf{16} (2008), 19--34.
\bibitem{C} A.~P. Calder{\'o}n, \emph{On an inverse boundary value problem}, Seminar on
  Numerical Analysis and its Applications to Continuum Physics, Soc. Brasileira
  de Matem{\'a}tica, R{\'i}o de Janeiro, 1980.
\bibitem{dOI} L.~d'Onofrio, T.~Iwaniec, \textit{{Notes on $p$-harmonic analysis}}, Contemp. Math. \textbf{370} (2005), 25--50.
\bibitem{HKM} J.~Heinonen, T.~Kilpel\"ainen, O.~Martio, Nonlinear potential theory of degenerate elliptic equations. Oxford Science Publications, Clarendon Press, Oxford, 1993.
\bibitem{HS} D.~Hervas, Z.~Sun, \textit{An inverse boundary value problem for quasilinear elliptic equations}, Comm. PDE \textbf{27} (2002), 2449--2490.
\bibitem{IUY} O.~Imanuvilov, G.~Uhlmann, M.~Yamamoto, \textit{{The Calder{\'o}n problem with partial data in two dimensions}}, J. Amer. Math. Soc. \textbf{23} (2010), 655--691.
\bibitem{I1} V.~Isakov, \textit{{On uniqueness in inverse problems for semilinear parabolic equations}}, Arch. Rat. Mech. Anal. \textbf{124} (1993), 1--12.
\bibitem{I2} V.~Isakov, \textit{{Uniqueness of recovery of some systems of semilinear partial differential equations}}, Inverse Problems \textbf{17} (2001), 607--618.
\bibitem{I3} V.~Isakov, \textit{{Uniqueness of recovery of some quasilinear partial differential equations}}, Comm. PDE \textbf{26} (2001), 1947--1973.
\bibitem{IN} V.~Isakov, A.~Nachman, \textit{{Global uniqueness for a two-dimensional semilinear elliptic inverse problem}}, Trans. Amer. Math. Soc. \textbf{347} (1995), 3375--3390.
\bibitem{IS} V.~Isakov, J.~Sylvester, \textit{{Global uniqueness for a semilinear elliptic inverse problem}}, Comm. Pure Appl. Math. \textbf{47} (1994), 1403--1410.
\bibitem{JK} D.~Jerison, C.~Kenig, \textit{{The inhomogeneous Dirichlet problem on Lipschitz domains}}, J. Funct. Anal. \textbf{130} (1995), 161--219.
\bibitem{KN} H.~Kang, G.~Nakamura, \textit{Identification of nonlinearity in a conductivity equation via the Dirichlet-to-Neumann map}, Inverse Problems \textbf{18} (2002), 1079--1088.
\bibitem{KY} H.~Kang, K.~Yun, \textit{{Boundary determination of conductivities and Riemannian metrics via local Dirichlet-to-Neumann operator}}, SIAM J. Math. Anal. \textbf{34} (2003), 719--735.
\bibitem{KSU} C.~E.~Kenig, J.~Sj\"ostrand, G.~Uhlmann,  \textit{The Calder\'on problem with partial data}, Ann. of Math. \textbf{165} (2007), 567--591.
\bibitem{KV} R.~Kohn, M.~Vogelius, \textit{{Determining conductivity by boundary measurements}}, Comm. Pure Appl. Math. \textbf{37} (1984), 289--298.
\bibitem{Kuijper} A.~Kuijper, \textit{{$p$-Laplacian driven image processing}}, Proc. International Conference of Image Processing, ICIP 2007 (San Antonio, Texas, USA), vol. V, 257--260.
%\bibitem{LeU} J.~Lee, G.~Uhlmann, \textit{Determining anisotropic real-analytic conductivities by boundary measurements}, Comm. Pure Appl. Math., \textbf{42} (1989), 1097--1112.
\bibitem{LW}
X.~Li, J.-N.~Wang, \emph{{Determination of viscosity in the stationary Navier-Stokes equations}}, J. Diff. Eq. \textbf{242} (2007), 24–-39.
\bibitem{L}
P.~Lindqvist, \emph{{Notes on the $p$-Laplace equation}}, Rep. Univ. Jyv\"askyl\"a Dept. Math. Stat. \textbf{102} (2006), 80 pp.
\bibitem{N}
A.~Nachman, \emph{Reconstructions from boundary measurements}, Ann. Math. \textbf{128} (1988), 531--576.
\bibitem{N_2D}
A.~Nachman, \emph{{Global uniqueness for a two-dimensional inverse boundary value problem}}, Ann. of Math. \textbf{143} (1996), 71--96.
\bibitem{NT}
G.~Nakamura, K.~Tanuma, \textit{{Formulas for reconstructing conductivity and its normal derivatives at the boundary from the localized Dirichlet to Neumann map}}, Proceedings of the international conference on inverse problems - recent developments in theories and numerics (ed. Y.-C.~Hon, M.~Yamamoto, J.~Cheng, J.-Y.~Lee) (2003), World Scientific, 192--201.
\bibitem{Sun_nonlinear} Z.~Sun, \textit{On a quasilinear boundary value problem}, Math. Z. \textbf{221} (1996), 293--305.
\bibitem{Sun2004} Z.~Sun, \textit{{Inverse boundary value problems for a class of semilinear elliptic equations}}, Adv. in Appl. Math. \textbf{32} (2004), 791--800.
\bibitem{Sun_survey} Z.~Sun, \textit{Conjectures in inverse boundary value problems for quasilinear elliptic equations}, Cubo \textbf{7} (2005), no. 3, 65--74.
\bibitem{Sun2010} Z.~Sun, \textit{{An inverse boundary-value problem for semilinear elliptic equations}}, Electron. J. Diff. Eq. \textbf{37} (2010), 1--5.
\bibitem{SuU} Z.~Sun, G.~Uhlmann, \textit{Inverse problems in quasilinear anisotropic media}, Amer. J. Math. \textbf{119} (1997), 771--797.
\bibitem{SU} J.~Sylvester, G.~Uhlmann, \textit{A global uniqueness theorem for an inverse boundary value problem}, Ann. of Math. \textbf{125} (1987), 153--169.
\bibitem{SU2} J.~Sylvester, G.~Uhlmann, \textit{Inverse boundary value problems at the boundary -- continuous dependence}, Comm. Pure Appl. Math. \textbf{41} (1988), 197--219.
\bibitem{U_IP} G.~Uhlmann, \textit{Electrical impedance tomography and Calder{\'o}n's problem}, Inverse Problems \textbf{25} (2009), 123011.
\bibitem{W} T.H.~Wolff, \textit{{Gap series constructions for the $p$-Laplacian}}. J. d'Analyse Math. \textbf{102} (2007), 371--394. Originally appeared as a preprint in 1984.

\end{thebibliography}
\end{document}